\documentclass[11pt]{article}
\usepackage{fontspec}
\usepackage[english,bidi=default]{babel}
\usepackage{microtype}
\usepackage{geometry}
\usepackage{setspace}
\usepackage{csquotes}
\usepackage[hidelinks]{hyperref}
\usepackage{bookmark}
\usepackage[
  backend=biber,
  style=numeric-comp,
  sorting=nyt,
  maxbibnames=99
]{biblatex}

\babelprovide[import]{persian}
\babelfont[persian]{rm}{Amiri-Regular.ttf}
\begin{document}

\title{\textbf{Mathematics in the Age of the Reproduction of Intelligence}}
\author{Ali Enayat \\
%EndAName
Department of Philosophy, Linguistics, and Theory of Science\\
University of Gothenburg}
\date{\today }
\maketitle

\begin{abstract}
Drawing on Walter Benjamin's account of technological reproducibility and
aura in relation to art, together with insights from Yehuda Rav and Ludwig
Wittgenstein, this speculative essay reflects on the impact of AI systems on
the practice of pure mathematics.
\end{abstract}

\bigskip

\section{Prologue}

This is not a scholarly paper on the philosophy or history of mathematics.
Rather, it is an essayistic attempt by a mathematician whose research focus
is on logic and set theory, and who has long been sensitized to the
philosophical, historical, and aesthetic aspects of his work, to give shape
to his reflections surrounding the question 
\enquote{How will AI transform
mathematics?}. I am addressing readers, especially mathematicians, with a
certain predilection for reflecting on humanistic aspects of higher
mathematics.

This essay arose from my recent realization that Walter Benjamin's analysis
in \emph{The Work of Art in the Age of Its Technological Reproducibility}
(which I first read around 2011), Rav's assessment in \emph{Why Do We Prove
Theorems?} (which I encountered around 2004), and Wittgenstein's views on
surveyability in \emph{Remarks on the Foundations of Mathematics} (which I
met in my first year in graduate school, in the early 1980s) collectively
offer a promising framework for the task at hand.\footnote{%
After developing this connection, I became aware of Stephanie Dick's earlier
essay \cite{dick2011}, which invokes Benjamin in its historical study of the
Automated Reasoning Assistant AURA and of human--machine collaboration in
mathematical proof. Dick's analysis, however, proceeds in a different
direction.}

The shape of the essay has also been influenced by insights gained from the
phenomenological view of mathematics. Thus a key concern here is with the
ways mathematical objects, proofs, and practices are experienced and made
intelligible. Among mathematicians, Gian-Carlo Rota developed an explicitly
phenomenological approach to mathematical experience, notably in his
discussions of mathematical truth, beauty, and proof \cite%
{rota1991truth,rota1997beauty,rota1997proof}.\footnote{%
For readers unfamiliar with phenomenology: the tradition was founded by
Edmund Husserl (1859--1938), who was originally trained as a mathematician,
studied in Berlin with Karl Weierstrass and Leopold Kronecker and received
his doctorate in Vienna in 1883. Phenomenology, in the broad sense relevant
here, concerns the structures of experience and the ways things appear to
consciousness. For an introduction to Husserl and phenomenology, see \cite%
{zahavi2025husserl}; on the importance of mathematics for Husserl's
philosophical development, see \cite{hartimo2021husserl}.} Against this
backdrop, I suggest two independent perspectives on the ways AI can affect
mathematical experience and practice. The Benjamin line emphasizes the
aesthetic and cultural components of mathematics; the Wittgenstein--Rav line
emphasizes the various functions of proofs in the development of mathematics
that go beyond their service as certificates of veracity.\footnote{%
After developing the framework of this essay, I became aware that Terence
Tao had raised closely related concerns in his July 2026 ICM public lecture, 
\emph{Mathematics in the Age of AI}, subsequently developed in a paper of
the same title \cite{tao2026icm,tao2026ai}. Tao's approach is different, but
his concern with what happens when human understanding and absorption are
overwhelmed by the production and verification of proofs overlaps with the
analysis in Section 5.}

In this essay I am principally concerned with the impact of AI on the \emph{%
practice of} \emph{pure mathematics}. Contrary to their stereotype as `cold
thinkers', pure mathematicians are guided not only by rigor, but also by
aesthetic criteria in choosing problems and formulating solutions---a
laborious process that demands considerable emotional resilience.\footnote{%
Donald R. Weidman drew attention to the emotional perils of mathematical
work in his brief 1965 letter to \emph{Science} \cite{weidman1965} (later
reprinted in the anthology \emph{Mathematics: People, Problems, Results} 
\cite[pp.~289--290]{campbellhiggins1984}). Weidman's letter prompted me to
`converse' with him in Persian in 2025 \cite{enayat2025perilsfa}; for an
English translation, see \cite{enayat2026perils}.} Despite having a
supremely rigorous and regimented framework, mathematics has a creative,
artistic core. Mathematicians often praise each other's accomplishments in
terms that carry unmistakable aesthetic judgments: `brilliant,' `beautiful,'
`elegant,' and `surprising.' A good piece of new mathematics should reveal a
striking relationship that makes us wiser. Like the best works of poetry, it
should unveil the deepest insights within the confines of a few strokes. In
the words of the nineteenth-century English mathematician James Joseph
Sylvester: 
\enquote{May not Music be described as the Mathematic of sense (sic),
Mathematic as Music of the reason? the soul of each the same!} This
perspective does not place mathematics within the classical domain of the
arts, but it makes Benjamin's analysis less remote from mathematical
practice than it may appear at first sight. What fundamentally separates the
two domains is that art's audience extends far beyond artists, whereas
high-level mathematical research is primarily appreciated by a specialized
cadre of fellow experts.

In a maneuver borrowed from classical martial arts, Baron von M\"{u}%
nchhausen, and G\"{o}del, the subject's own strength is reflected back on
itself in the preparation of this essay. The first draft was developed
through an extended process of author-directed, AI-assisted drafting,
criticism, and revision. I conceived the project, determined its
intellectual direction, and supplied the principal arguments and critical
interventions. ChatGPT, in turn, assisted in developing and testing
formulations, reorganizing the argument, identifying relevant literature,
checking bibliographic details, and preparing and troubleshooting the LaTeX
source.\footnote{%
The version of ChatGPT used in preparing this essay was made available to me
through the University of Gothenburg.} However, in the preparation of
subsequent drafts, I worked on the text anew with minimal AI assistance in
order to compress several parts and incorporate more of my personal views,
experiences, and prose style. I also benefited from substantive comments and
suggestions provided by Amy Enayat, and from collegial support of Sergio
Fajardo and Siavash Shahshahani.

\section{Mathematics through History}

\label{mathematics-as-a-historical-practice}

Before turning to the main ideas of the paper, let us briefly review how
mathematical practice has evolved through history, and how pure mathematics,
with its supreme attention to internal abstract considerations and logical
harmony, took shape in its current form, a form now facing the
transformative currents of AI technology.

Mathematics is often described as the exemplar of timelessness and
independence from human contingencies. The Pythagorean theorem, for example,
appears to be independent of the foibles of language, politics, technology,
and fashion. Yet even its very name indicates a subtle tension that is
addressed in this essay: the proposition's validity is independent of
Pythagoras, while its cultural identity is not. In a less obvious way, the
very word \textquotedblleft algorithm\textquotedblright\ preserves, in
Latinized form, the name of al-Khwarizmi, the Persian polymath whose works
were instrumental in transmitting both algebraic methods and Hindu-Arabic
numeration.

The earliest surviving records reveal that mathematics arose as an
inseparable component of certain practices. In Mesopotamia and Egypt it was
connected with accounting, taxation, measurement, construction, and
astronomy; Chinese and Indian traditions similarly developed sophisticated
mathematical techniques in response to problems of administration, trade,
calendrical calculation, and astronomy \cite{katz2009,joseph2011}. Yet these
practical concerns were never entirely separable from questions of cosmic
order and meaning. For example, in many ancient societies, our current
distinction between astronomy as a hard science and astrology as a source of
esoteric knowledge would have been meaningless.

The subsequent history of mathematics is full of movements between
externally prompted problems and internally generated developments. Algebra
in the Islamic world was developed in a context that included inheritance,
trade, surveying, and astronomical calculation, yet it became a theory with
its own internal questions. Trigonometry developed through astronomy and
navigation but became a mathematical domain in its own right. Probability
emerged partly from games of chance, insurance, and demographic concerns but
became central to statistics, physics, economics, and pure mathematics.
Calculus arose in intimate connection with problems of motion, tangency,
areas, and celestial mechanics, only to become an immense theoretical
edifice whose concepts now reach far beyond those origins.

Over time, internally generated mathematical questions increasingly acquired
an autonomous life, developing into the study of abstract structure and
logical rigor, with proofs standing out as the principal currency of
correctness. Meanwhile, the empirical tools used for land surveying,
accounting, and celestial tracking matured into the systematic modeling of
physical phenomena. The division between \textquotedblleft applied
mathematics\textquotedblright\ and \textquotedblleft pure
mathematics\textquotedblright\ did not arise until the 18th and 19th
centuries to distinguish between them. Nowadays, however, this distinction
has become increasingly porous. Fields once pursued largely for internal
theoretical reasons---number theory, abstract algebra, topology, and logic
among them---have found major applications in cryptography, computing,
physics, and other sciences, making the boundary between \textquotedblleft
pure\textquotedblright\ and \textquotedblleft applied\textquotedblright\
mathematics far less clear-cut than it once appeared. Indeed, von Neumann
spoke of a \textquotedblleft quite peculiar duplicity\textquotedblright\ in
the nature of mathematics, i.e., a Janus-like feature, where one face finds
its sources of inspiration in experimental and applied contexts, while the
other one pursues purely internal, abstract considerations \cite%
{vonneumann1947mathematician}.

\section{Proof, Formalization, and Computation}

\label{proof-formalization-and-computation}

This section offers a `telegraphic' review of how proof-based mathematics
came to be formalized and mechanized. The later sections do not explicitly
depend on this history, but it helps place recent advances in AI within a
wider context in which formal logic and computation became increasingly
prominent.

Greek mathematics introduced a model of deductive proof whose later
influence is difficult to exaggerate. Euclid's \emph{Elements} did not
merely collect useful results; it organized geometry as a deductive
structure in which propositions followed from definitions, postulates, and
earlier propositions. Although actual mathematical practice has always been
more varied than the Euclidean ideal suggests, proof became one of the
defining features through which mathematics distinguished itself from other
forms of knowledge. A proof of a theorem shows that a certain proposition is
not merely supported by observation but that, under the relevant
assumptions, the conclusion could not be otherwise.

The modern project of formalization can reasonably be said to have entered a
decisive new phase with Frege's \emph{Begriffsschrift} (1879), which
introduced a formal language capable of representing quantification and
complex inferential structure, while his broader logicist program sought to
ground arithmetic in logic \cite{frege1967}. Whitehead and Russell's \emph{%
Principia Mathematica} (1910--1913) then attempted, on a monumental scale,
to reconstruct extensive portions of mathematics from explicitly stated
logical principles \cite{whiteheadrussell1910}. Whatever the fate of
logicism, these works helped establish a decisive possibility: mathematical
reasoning of various sorts could be symbolically represented within formal
systems. Hilbert took the project of formalization in a different direction.
His proof-theoretic program aimed to formalize mathematical theories and
then study those formal systems metamathematically, with the hope of
securing their consistency by suitably restricted, finitary means \cite%
{hilbert1983}. Unlike logicism, the point was not to identify mathematics
with logic; it was to make formal proofs themselves objects of mathematical
investigation. This sharpened a distinction that remains important today: a
derivation may become increasingly explicit as a symbolic object without
necessarily becoming more illuminating to a mathematician.

G\"{o}del's incompleteness theorems of 1931 transformed this foundational
landscape. For any consistent, effectively axiomatized formal system capable
of modest arithmetical strength, there are arithmetical statements that
cannot be proved within the system; and under a slightly stronger assumption
of arithmetical strength such a system cannot even prove the arithmetical
sentence that expresses the system's consistency. These results did not show
that mathematics was unreliable, nor did they make formal methods futile.
But they demonstrated that the hope of capturing sufficiently rich
mathematics within a complete and internally self-certifying formal
framework faced logical limitations \cite{godel1931}. At about the same
time, another question was taking shape: what precisely is an effective (or
mechanical) procedure? The independent work of several researchers,
including Post, Church, and Turing, supplied decisive mathematical
formulations of effective calculability. Turing's 1936 account described an
idealized machine that manipulates symbols according to finite rules \cite%
{church1936,turing1937}. The Turing machine was not conceived as an
engineering design for an electronic computer, but as an abstract device for
giving mathematical precision to the intuitive notion of mechanical
calculation. This episode established an entanglement that would deepen over
the next century: logic contributed to computer science; computers
transformed mathematical practice; formal logic later became executable in
proof assistants; and contemporary AI systems based on neural-network
technology now use mathematics both as a domain of application and as a
training ground for general reasoning. \emph{This historical trajectory
reveals key shifts in the role of formal logic and computation in
mathematics---from explicit formal representation, to interactive proof
assistants }(\emph{which blend formal logic with computational execution})%
\emph{, and ultimately to neural-network-generated proofs verified by humans
or automated proof checkers.}

\section{Benjamin: The Dispersal of Aura}

\label{benjamin-reproducibility-aura-and-phenomenology}

To take the first steps in probing the possible impact of the emergence of
machine-generated mathematical reasoning on the \emph{lived mathematical
experience}, I have found it useful to turn to a thinker who, at first
sight, seems far removed from mathematics and computation. Walter Benjamin
(1892--1940) was a German philosopher, literary critic, and cultural
theorist who wrote about aesthetics, language, history, technology,
architecture, photography, and modern urban life, often interweaving all of
them. He was intellectually close to the founders of the Frankfurt School,
though his position was always idiosyncratic. One of his most well-known
works is \enquote{The
Work of Art in the Age of Mechanical Reproduction} \cite{benjamin2008work}.
He wrote versions of his influential essay after being forced into exile
following the Nazi seizure of power in the mid-1930s, a period during which
photography and cinema became pervasive, and witnessed the emergence of mass
politics, propaganda, and fascism. He died in 1940 while attempting to
escape Nazi-occupied Europe. The essay probes the ways in which
technological innovations in the reproduction and circulation of artworks
transform their value, meaning, and modes of reception. The transformation
itself was set in motion decades prior to Benjamin's essay; by contrast, we
appear to be at the threshold of the technological reproduction of
intelligence.

Benjamin's analysis begins from a seemingly straightforward distinction: a
traditional artwork such as a painting or sculpture has a \textit{unique
material existence}. It is understood as this object, in this place, with
this physical history. Copies may represent it, but they are mere pointers
to it. For Benjamin, that specificity is related to the work's \emph{aura}:
a difficult and deliberately suggestive term that evokes uniqueness,
distance, authenticity, historical embeddedness, and the authority of
tradition. The aura of a work is not circumscribed by its fame or prestige;
rather, it concerns the kind of relation we have to an object that can be at
once near and remote, not geometrically, but semantically. Benjamin argues that technological reproduction disrupts this relation: as reproductions multiply and circulate independently of the original, the work's singular ‘here and now’ loses its privileged status. In his own words, \enquote{what withers in the age of the technological reproducibility
of the work of art is the latter's aura} \cite[p.~22]{benjamin2008work}. \footnote{
There is a phenomenological texture to Benjamin's notion of aura, since it
concerns the manner in which objects are experienced or appear to us.
Benjamin was not a Husserlian phenomenologist, but aura is also a
description of how a relation is lived. As argued by Miriam Bratu Hansen,
aura is not simply a property of an object, but a structure of experience 
\cite{hansen2008}.}

Benjamin's point is therefore far stronger than the claim that reproduction
lowers quality or produces cheap copies. Technology changes the very
experience of reception and, with it, the social function of the work. Thus 
\emph{cult value}, grounded in ritual, singular presence, and meaning,
undergoes a metamorphosis into \emph{exhibition value}, grounded in
circulation, display, and consumption. It is exemplified in the contrast
between a prolonged period of reflection in a monastery and the distracted
participation in a public ceremony. Benjamin juxtaposes this kind of
distraction with what one English translation calls \emph{contemplative
immersion}, during which the beholder is absorbed by the work; in distracted
reception, by contrast, the work is absorbed into the habitual activity of
its audience \cite[pp.~39--40]{benjamin2008work}. The former experience is
immediately familiar to mathematicians: prolonged immersion in an abstract
world, during which the surrounding environment recedes and the mathematical
object acquires a peculiar experiential nearness. To produce first-rate
results in pure mathematics, uninterrupted periods of time for study and
reflection, along with lots of paper (and a wastebasket for many false roads
traveled), have typically been of more value than an army of graduate
students, expensive laboratories or equipment.

In a later essay \emph{On Some Motifs in Baudelaire}, Benjamin adds the
striking feature of \emph{returning the gaze} to the auratic experience:
``To perceive the aura of an object we look at means to invest it with the
ability to look at us in return''. \cite[p.~338]{benjamin2003baudelaire}.
The pertinent point is that perception is no longer described as a one-way
act, but one that involves reciprocity: the contemplated object seems to
resist complete appropriation and to address the perceiver in return.%
\footnote{%
Arianne Conty interprets this dimension as approaching a theory of
intentionality organized around the reciprocal gaze; see \cite{conty2013}.}
There is a suggestive analogue in the language mathematicians often use to
describe their own experience: a complex problem or a theorem can seem to
resist, or to \textquotedblleft push back\textquotedblright\ against
attempts to force it into a preconceived form. In other words, it may appear
that it is imposing demands of its own, contrary to the mathematician's
expectations. There are many famous examples of this phenomenon throughout
the history of mathematics, typified in the following exchange. In answer to
the question \textquotedblleft What research problems and areas are you
likely to explore in the future?\textquotedblright, Maryam Mirzakhani \cite%
{Mirzakhani2008} answered \textquotedblleft It's hard to predict. But I
would prefer to follow the problems I start with wherever they lead
me.\textquotedblright

Prima facie, mathematics seems to shun even the very idea of aura. A proof
is supposed to displace personal authority: Euclid's lemma is not true
because Euclid was a wise sage, and Noether's theorem does not acquire
additional validity from Noether's pedigree. Ideally, anyone who has
mastered the relevant competence can reconstruct the argument. So while
mathematical truth projects an image devoid of human traces, there is an
aura-like texture to the way it is understood and remembered. Mathematical
culture continually restores history, personality, and genealogy to truths
whose validity does not require them. The names are everywhere: Euclidean
geometry, Euler characteristic, Hilbert space, Minkowski geometry, Noether's
theorem, G\"{o}del numbering, and the Poincar\'{e} disk. Such eponyms are
not always historically accurate records of priority, and mathematicians
often debate whether they distribute credit fairly, but precisely for that
reason they are revelatory. They do cultural and phenomenological work in
that they attach an otherwise reproducible object to a singular lineage of
discovery, a style of thought, a school, or a remembered person.
AI-generated proofs disturb this genealogy---and thereby the associated
aura---by dispersing the roots of mathematical genesis across
mathematicians, models, formal libraries, and institutions. \emph{Thus,
theorems generated by AI may possess a perfectly legitimate proof while
lacking a decipherable history and a recognizable `persona'.}

\section{The Function of Proof: Wittgenstein and Rav}

\label{proof-as-knowledge-surveyability-rav-and-formal-verification}

In \emph{Remarks on the Foundations of Mathematics}, Wittgenstein made the
surveyability (\emph{\"{U}bersichtlichkeit}) of mathematical proof a central
concern, closely linking it to reproducibility and to our ability to
recognize whether the same proof has been reproduced. The point is not that
a proof should be short, easy, or grasped at once; rather, it must be
presented in a form that can function as a recognizable and reproducible
object within human mathematical practice \cite%
{wittgenstein1978,muehlhoelzer2005}.

The 1976 computer-assisted proof of the Four Color Theorem by Kenneth Appel
and Wolfgang Haken put the notion of surveyability to the test. Its
machine-checked cases could be rerun, but no individual mathematician could
realistically traverse them in the ordinary manner of reading a proof from
beginning to end \cite{appel1977a,appel1977b}. Thomas Tymoczko's distinctive
analysis treated the episode not merely as an unusually long proof but as a
challenge to commonly held assumptions about the epistemological character
of mathematical proof \cite{tymoczko1979}. The controversy made explicit a
distinction that remains implicit in ordinary mathematical practice: a proof
can certify that a conclusion follows from accepted assumptions without
perspicuously explaining why it holds. A particularly useful way to
understand the various functions of mathematical proofs was proposed by
Yehuda Rav in his 1999 essay \textquotedblleft Why Do We Prove
Theorems?\textquotedblright\ \cite{rav1999}. Rav distinguishes ordinary
mathematical proofs from formal derivations and argues that the epistemic
content of a proof goes far beyond the proposition summarized in the
theorem. Proofs carry methods, constructions, strategies, conceptual
distinctions, representations, and connections to other problems. On this
view, a theorem is often the synoptic public face of a multifaceted process
whose mathematical significance cannot be reduced to the truth of its
concluding sentence. To understand why mathematicians prove theorems, one
must therefore ask what is learned in and through the proof, not merely that
the proposition has been certified.

Rav dramatizes the point with an imaginary universal decision machine
marketed under the name PYTHIAGORA. Rav imagines, contrary to the well-known
limits of decidability, that a mathematician could type any appropriately
formulated problem into an oracular computer and immediately receive the
answer `true' or `false'.\ If the sole purpose of proof were to ascertain
truth, PYTHIAGORA would make proving largely redundant: mathematicians could
devote themselves to conjecture and let the machine settle correctness.
Rav's point is that this conclusion rests on too narrow a view of what
proofs do. The proof is not disposable scaffolding around the theorem. It is
one of the principal junctures at which mathematical knowledge is
articulated and enlarged.

Rav's account explains familiar features of actual mathematical culture.
Mathematicians routinely seek new proofs of theorems already known to be
true. A second proof can be more valuable than the first because it reveals
a different structure, introduces a technique that applies elsewhere,
removes an unnecessary hypothesis, connects two fields, suggests a
generalization, or makes a previously mysterious result intelligible. It is
one thing to prove a theorem, and a whole other one to find the `right'
proof for it. In other words, a proof may very well correctly establish a
result while being judged lacking if it offers little guidance about why the
result holds or what to look for next. The value of proof therefore does not
reside solely in its function as a certificate, but also in its capacity to
lead to human understanding and enlightenment.

In this day and age, we find ourselves---unexpectedly---in a situation in
which Rav's imaginary oracle appears rather technologically modest.
Contemporary AI systems do not merely decide whether certain propositions
are true; they can propose intermediate lemmas, search proof spaces,
translate informal arguments into formal ones, generate alternative
explanations, suggest analogies, and in some settings construct substantial
portions of proofs. The relevant question is no longer what mathematics
would become if proofs were relegated to an oracle. It is what mathematics
becomes if proofs themselves, including parts of their generative structure,
can be supplied at the touch of a button, with no group of humans to pose
the question ``how did you come up with that?''.

Machine-generated proofs may introduce constructions and strategies from
which human mathematicians genuinely learn---methods that are reused,
generalized, and taught. However, the outcome that appears likely is a
proliferation of AI-generated esoteric proofs, a significant number of which
will be proved correct by proof assistants. Certified proofs would become
increasingly abundant while human understanding, transferable methods, and
conceptual organization would woefully lag behind. Formal verification makes
the underlying tension especially sharp: correctness can be certified
without guaranteeing mathematical understanding.\footnote{%
For complementary discussions, see Venkatesh \cite%
{venkatesh2024automation,venkatesh2025human}. The former examines how
automation may alter the processes by which mathematical value and consensus
are formed; the latter emphasizes the conceptual and communicative
structures through which mathematics becomes humanly intelligible.}

This suggests updated Wittgensteinian questions: \textit{Will surveyability
disappear through the diffusion of the proof to an intractable set of
generative sources? Can we reasonably hope that reproducibility and checking
only partly migrate to machines while human assessment remains the final
arbiter of mathematical understanding and significance}?

\section{Concluding Thoughts}

The framework developed here is intended to sketch some of the fault lines
along which the mathematical experience may be revolutionized. The
transformation we are witnessing is of colossal proportions, and its
eventual form remains profoundly unpredictable; as in the attempt to
understand a large-scale geological disruption, any conceptual framework can
at best probe its contours. It is as if one tried to anticipate the
long-term effects on the practice of calligraphy on the eve of Gutenberg's
invention of printing.

Due to turbulent socio-political and economic forces, and despite the best
efforts within the mathematical community, the transition is unlikely to be
orderly or predictable. Mathematics has entered a tumultuous period in which
established relations among truth, credit, provenance, understanding,
judgment, and intellectual responsibility are repeatedly unsettled before
new conventions can take shape. Some old hands---myself squarely among
them---will strive to preserve domains in which human discovery, beauty,
surveyability, and understanding are non-negotiable desiderata. Others will
participate fully in distributed forms of human--machine mathematical
collaboration, while others will move uneasily between these extremes.

In the short run, mathematicians will increasingly employ AI---akin to
children in a candy store---to unravel many longstanding enigmas, along with
building unseen bridges between seemingly disconnected mathematical
territories. This will be accomplished, of course, by drawing on the vast existing corpus of mathematical knowledge amassed by human ingenuity, which serves both as training material for language models and as a resource available to AI agents. But what happens when AI advances beyond the frontier of humanly digested mathematics, and increasingly builds upon AI-generated `knowledge'? This informs a simile that I recently shared with a colleague:
\textquotedblleft At the moment, I am more thrilled than scared, as if I
were swimming in the ocean, and a massive wave has suddenly picked me up and
is taking me to unbelievable heights. I am very much enjoying the ride, but
I know all too well that I cannot fathom the large-scale impact of such a
colossal wave when it crashes\textquotedblright .

In the long run, we may find ourselves in a situation reminiscent of Stanis\l{}aw Lem's novel \emph{Solaris} \cite{lem2011solaris}, which depicts the struggle of scientists—dubbed \emph{Solarists}—with understanding a planet covered by a mysterious sentient ocean. In the novel, \emph{Solaristics} has evolved over decades into an entire academic discipline, complete with its own history, competing theories, and specialized departments—all dedicated to trying (and largely failing) to comprehend or communicate with the alien planet. Benjamin's reflections on photography offer a counterpoint to Lem's nightmarish image. Technology, he points out,
can disclose aspects of reality inaccessible to unaided perception:
photography, through devices such as enlargement and slow motion, reveals
what he calls the \enquote{optical unconscious} \cite[pp.~511--512]%
{benjamin1931photography}. Similarly, the AI lens may enable mathematicians
to perceive hitherto inaccessible structures and analogies, and perhaps to
develop altogether new modes of mathematical cognition.

\emph{It is unclear what the future holds for younger mathematicians,
especially because---as I suspect---advances in AI will ultimately lead to
the appropriation of pure mathematics by the market forces of applied
mathematics. That is, we may see the further dominance of profitable
veracity over edifying beauty; or, if you will, techne over ars. We may well
be entering a `neo-Babylonian era', where mathematics is once again
inseparable from utility.}

\printbibliography[title={References}]

\end{document}